\documentclass[amssymb,aps,showpacs,showkeys,showpacs,notitlepage]{revtex4-1}
\usepackage{graphicx}
\usepackage{amsmath}
\usepackage[pdfauthor={Richard J. Mathar},pdfkeywords={Errata, Exchange Integral, Two-dimensional Fermi Gas},pdftoolbar=false,colorlinks=true]{hyperref}
\DeclareMathOperator{\csch}{csch}

\begin{document}

\title[PV of a quadruple integral]{Erratum to ``Solutions problem 89-2: On the principal value of a quadruple integral'', SIAM Rev. 32 (1990) 143.}

\author{Richard J. Mathar}
\homepage{http://www.mpia.de/~mathar}
\affiliation{Max-Planck Institute of Astronomy, K\"onigstuhl 17, 69117 Heidelberg, Germany}

\date{\today}

\begin{abstract}
W.\ B.\ Jordan's conclusion that the quadruple principal value integral in problem 89-2 vanishes
does not hold. The error sneaks in through a contribution of a subintegral which impedes some
sign symmetry with respect to the master parameter (the Fermi radius) and which was
overlooked in the published solution. In summary, the
original problem of solving the quadruple integral remains
unsolved.
\end{abstract}

\maketitle
\section{Statement of the Task}
\begin{figure}[hbt]
\includegraphics[scale=0.45]{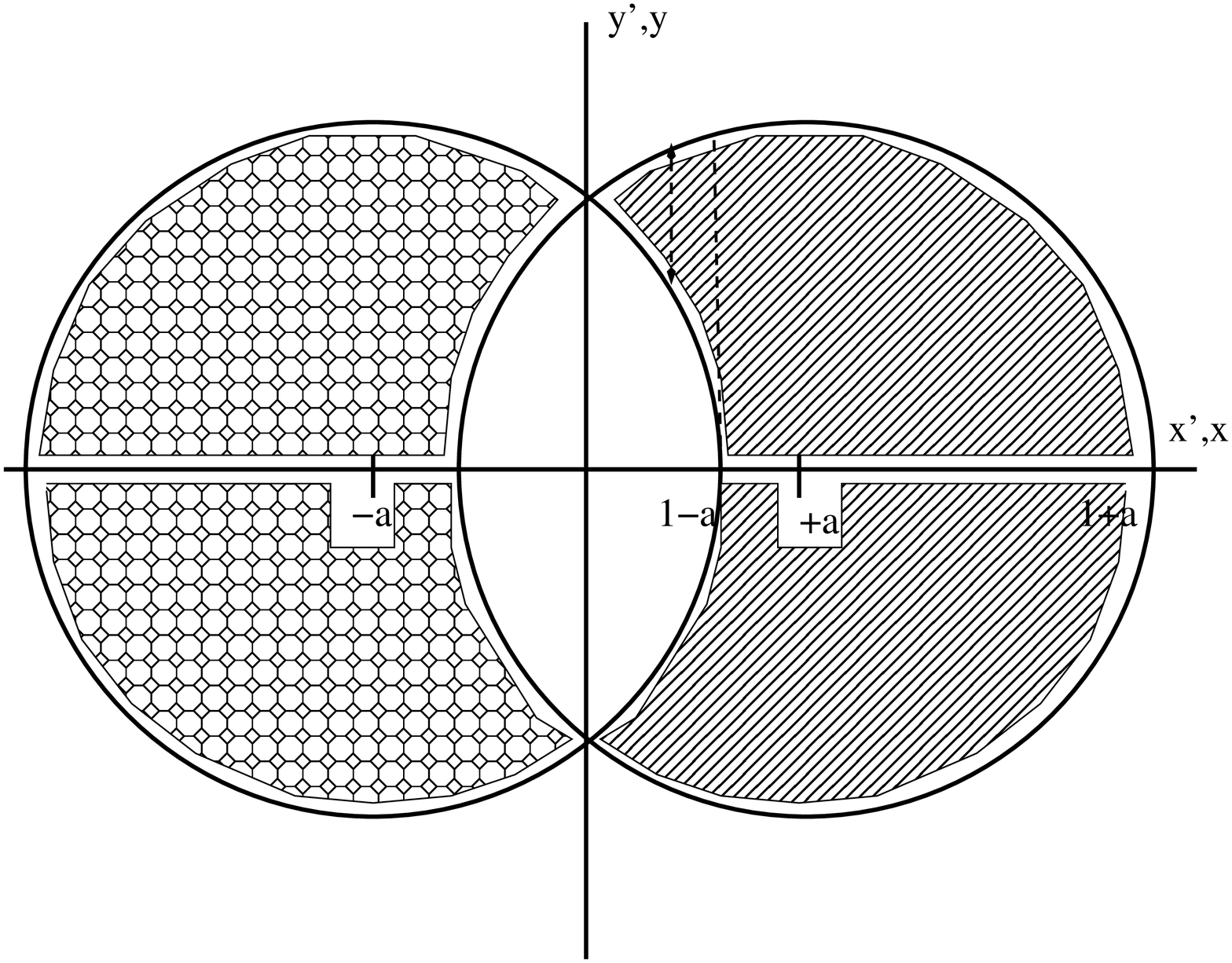}
\caption[]{The two moon-shaped regions of integration for $x'$ and $y'$
and for $x$ and $y$, one filled with bubbles, the other with hatching depending on the sign.
The circles intersect the vertical axis at $y,y'=\pm \sqrt{1-a^2}$.
}
\label{fig.circ}
\end{figure}

The manuscript is concerned with the evaluation of the principal value of the quadruple integral
\cite{GlasserSIAMR31}
\begin{equation}
F(a) = P\int_{-\infty}^\infty \frac{dx}{x^2}
\int_{-\infty}^\infty \frac{dx'}{x'}
\int_{-1}^1 dy
\int_{-1}^1 dy'
\frac{x-x'}{(x-x')^2+(y+y')^2}\Delta(x,y)\Delta(x',y')
.
\label{eq.Fdef}
\end{equation}
The function $\Delta(x,y)$ is defined to be $+1$ in the moon-shaped region
inside the unit circle centered at $x=a$, and it is $-1$ in the moon-shaped
mirror region
inside the unit circle centered at $x=-a$.
In the infinite exterior region and in the lens-shaped region
around the center, where the two circles overlap, $\Delta(x,y)$ is zero.
The two connected regions that contribute to the integral are illustrated in Figure \ref{fig.circ}. 
They touch each other at two common points on the vertical axis.

The two circles represent Fermi disks in the application
to solid state physics
\cite{MaldagueSSC26,RajagopalPRB15,GeldartCJP48_155};
this is not relevant to the further calculation.

The effect of switching the sign of the parameter $a$  is
\begin{equation}
\Delta(x,y) \to -\Delta(x,y);\quad
\Delta(x',y') \to -\Delta(x',y').
\end{equation}
So the product $\Delta(x,y)\Delta(x',y')$ is invariant towards changing the sign of $a$,
and
\begin{equation}
F(a)=F(-a)
.
\end{equation}
In that sense one only needs to consider $a\ge 0$
for the rest of the calculation.

\section{Criticism of Jordan's conclusion}

Jordan's calculation \cite{JordanSIAMR32}
argues that the double inner integral $I(x,y)$ of $F$, an
integral over the full right circle of $x'$ and $y'$ in Figure \ref{fig.circ},
is an even function of $a$. This would cause the entire integral
to vanish at the time when the sum over both circles
is involved, because $\Delta(x',y')$ is an odd function of $a$.

The error appears implicitly in the step
\begin{eqnarray*}
I_2/\pi &=& \int_0^a (A-Be^{-u})\csch u dr + \int_a^1 B dr
\\
&=& \int_0^a (A\csch u -B \coth u)dr +\int_0^1 Bdr,
\end{eqnarray*}
in the 4-th but last equation on page 144,
which splits $\int_0^a B dr$ off the
left integral and unites it with the right integral.
Although formally correct,
this step \emph{only} applies to the cases where $a\ge 0$. If $a$ is \emph{negative},
the upper limit of the first and the
lower limit of the second integral in the first of these two lines 
must be clamped to zero; this is basically a consequence of the role
of $r$ as a radial circular coordinate which cannot become negative.

By an equivalent reasoning, the step
is not correct if $a>1$, because then the upper limit on $r$ is 1
and the second integral $\int_a^1 Bdr$ should not contribute at all.

We see that for negative $a$ the first
integral $\int_0^a (A-Be^{-u})\csch u dr$ should not contribute at all
and the second be changed to $\int_0^1 Bdr$, so
the contribution of $\log[(p+2a)/p']$ to $I_2/\pi$ in the last equation
of page 144 must be dropped for negative $a$.

As a consequence, Jordan's final equation
\[
I=\frac{\pi}{2}\log(2x)^4[(x+iy)^2+1-a^2][(x-iy)^2+1-a^2]
\]
is invalid whenever $a$ is negative.
Although $I$ in that form is an even function of $a$, its validity is restricted
to $a>0$---and for the congruential reason concerning the upper limit also
to $a<1$. The cancellation claimed by Jordan when the values
at positive and negative $a$ are subtracted while calculating
the double integral over $x'$ and $y'$ is simply inhibited
because the values of $I_2$ at negative $a$ are not those obtained
by symmetric (even) extrapolation.

Glasser's quest \cite{GlasserSIAMR31} of obtaining values of the quadruple integral---perhaps
not analytically but merely with satisfactory numerical methods---remains
unanswered so far.

\bibliography{all}

\end{document}